\documentclass[11pt, a4paper]{article}
\usepackage[cp1251]{inputenc}
\usepackage[russian]{babel}
\usepackage{amsmath} \usepackage{euscript}
\usepackage{mymatrix}%\let\mymatrixII\mymatrix
\usepackage{longtable}
\usepackage{mathrsfs}
\usepackage{amssymb}
\begin{document}

\newcounter{bnomer} \newcounter{snomer}
\newcounter{bsnomer}
\setcounter{bnomer}{0}
\renewcommand{\thesnomer}{\thebnomer.\arabic{snomer}}
\renewcommand{\thebsnomer}{\thebnomer.\arabic{bsnomer}}
\renewcommand{\refname}{\begin{center}\large{\textbf{References}}\end{center}}

\setcounter{MaxMatrixCols}{14}

\newcommand{\sect}[1]{%
\setcounter{snomer}{0}\setcounter{bsnomer}{0}
\refstepcounter{bnomer}
\par\bigskip\begin{center}\large{\textbf{\arabic{bnomer}. {#1}}}\end{center}}
\newcommand{\sst}{%
\refstepcounter{bsnomer}
\par\bigskip\textbf{\arabic{bnomer}.\arabic{bsnomer}. }}
\newcommand{\defi}[1]{%
\refstepcounter{snomer}
\par\medskip\textbf{Definition \arabic{bnomer}.\arabic{snomer}. }{#1}\par\medskip}
\newcommand{\theo}[2]{%
\refstepcounter{snomer}
\par\textbf{Theorem \arabic{bnomer}.\arabic{snomer}. }{#2} {\emph{#1}}\hspace{\fill}$\square$\par}
\newcommand{\mtheop}[2]{%
\refstepcounter{snomer}
\par\textbf{Theorem \arabic{bnomer}.\arabic{snomer}. }{\emph{#1}}
\par\textsc{Proof}. {#2}\hspace{\fill}$\square$\par}
\newcommand{\mcorop}[2]{%
\refstepcounter{snomer}
\par\textbf{Corollary \arabic{bnomer}.\arabic{snomer}. }{\emph{#1}}
\par\textsc{Proof}. {#2}\hspace{\fill}$\square$\par}
\newcommand{\mtheo}[1]{%
\refstepcounter{snomer}
\par\medskip\textbf{Theorem \arabic{bnomer}.\arabic{snomer}. }{\emph{#1}}\par\medskip}
\newcommand{\theoc}[2]{%
\refstepcounter{snomer}
\par\medskip\textbf{Theorem \arabic{bnomer}.\arabic{snomer}. }{#1} {\emph{#2}}\par\medskip}
\newcommand{\mlemm}[1]{%
\refstepcounter{snomer}
\par\medskip\textbf{Lemma \arabic{bnomer}.\arabic{snomer}. }{\emph{#1}}\par\medskip}
\newcommand{\mprop}[1]{%
\refstepcounter{snomer}
\par\medskip\textbf{Proposition \arabic{bnomer}.\arabic{snomer}. }{\emph{#1}}\par\medskip}
\newcommand{\theobp}[2]{%
\refstepcounter{snomer}
\par\textbf{Theorem \arabic{bnomer}.\arabic{snomer}. }{#2} {\emph{#1}}\par}
\newcommand{\theop}[2]{%
\refstepcounter{snomer}
\par\textbf{Theorem \arabic{bnomer}.\arabic{snomer}. }{\emph{#1}}
\par\textsc{Proof}. {#2}\hspace{\fill}$\square$\par}
\newcommand{\theosp}[2]{%
\refstepcounter{snomer}
\par\textbf{Theorem \arabic{bnomer}.\arabic{snomer}. }{\emph{#1}}
\par\textbf{Схема доказательства}. {#2}\hspace{\fill}$\square$\par}
\newcommand{\exam}[1]{%
\refstepcounter{snomer}
\par\medskip\textbf{Example \arabic{bnomer}.\arabic{snomer}. }{#1}\par\medskip}
\newcommand{\deno}[1]{%
\refstepcounter{snomer}
\par\textbf{Definition \arabic{bnomer}.\arabic{snomer}. }{#1}\par}
\newcommand{\post}[1]{%
\refstepcounter{snomer}
\par\textbf{Proposition \arabic{bnomer}.\arabic{snomer}. }{\emph{#1}}\hspace{\fill}$\square$\par}
\newcommand{\postp}[2]{%
\refstepcounter{snomer}
\par\medskip\textbf{Proposition \arabic{bnomer}.\arabic{snomer}. }{\emph{#1}}%
\ifhmode\par\fi\textsc{Proof}. {#2}\hspace{\fill}$\square$\par\medskip}
\newcommand{\lemm}[1]{%
\refstepcounter{snomer}
\par\textbf{Lemma \arabic{bnomer}.\arabic{snomer}. }{\emph{#1}}\hspace{\fill}$\square$\par}
\newcommand{\lemmp}[2]{%
\refstepcounter{snomer}
\par\medskip\textbf{Lemma \arabic{bnomer}.\arabic{snomer}. }{\emph{#1}}
\par\textsc{Proof}. {#2}\hspace{\fill}$\square$\par\medskip}
\newcommand{\coro}[1]{%
\refstepcounter{snomer}
\par\textbf{Corollary \arabic{bnomer}.\arabic{snomer}. }{\emph{#1}}\hspace{\fill}$\square$\par}
\newcommand{\mcoro}[1]{%
\refstepcounter{snomer}
\par\textbf{Corollary \arabic{bnomer}.\arabic{snomer}. }{\emph{#1}}\par\medskip}
\newcommand{\corop}[2]{%
\refstepcounter{snomer}
\par\textbf{Corollary \arabic{bnomer}.\arabic{snomer}. }{\emph{#1}}
\par\textsc{Proof}. {#2}\hspace{\fill}$\square$\par}
\newcommand{\nota}[1]{%
\refstepcounter{snomer}
\par\medskip\textbf{Замечание \arabic{bnomer}.\arabic{snomer}. }{#1}\par\medskip}
\newcommand{\propp}[2]{%
\refstepcounter{snomer}
\par\medskip\textbf{Proposition \arabic{bnomer}.\arabic{snomer}. }{\emph{#1}}
\par\textsc{Proof}. {#2}\hspace{\fill}$\square$\par\medskip}
\newcommand{\hypo}[1]{%
\refstepcounter{snomer}
\par\medskip\textbf{Conjecture \arabic{bnomer}.\arabic{snomer}. }{\emph{#1}}\par\medskip}
\newcommand{\prop}[1]{%
\refstepcounter{snomer}
\par\textbf{Proposition \arabic{bnomer}.\arabic{snomer}. }{\emph{#1}}\hspace{\fill}$\square$\par}

\newcommand{\Ind}[3]{%
\mathrm{Ind}_{#1}^{#2}{#3}}
\newcommand{\Res}[3]{%
\mathrm{Res}_{#1}^{#2}{#3}}
\newcommand{\epsi}{\epsilon}
\newcommand{\tri}{\triangleleft}
\newcommand{\Supp}[1]{%
\mathrm{Supp}(#1)}

\newcommand{\lee}{\leqslant}
\newcommand{\gee}{\geqslant}
\newcommand{\reg}{\mathrm{reg}}
\newcommand{\empr}[2]{[-{#1},{#1}]\times[-{#2},{#2}]}
\newcommand{\sreg}{\mathrm{sreg}}
\newcommand{\codim}{\mathrm{codim}\,}
\newcommand{\chara}{\mathrm{char}\,}
\newcommand{\rk}{\mathrm{rk}\,}
\newcommand{\chr}{\mathrm{ch}\,}
\newcommand{\id}{\mathrm{id}}
\newcommand{\Ad}{\mathrm{Ad}}
\newcommand{\col}{\mathrm{col}}
\newcommand{\row}{\mathrm{row}}
\newcommand{\low}{\mathrm{low}}
\newcommand{\pho}{\hphantom{\quad}\vphantom{\mid}}
\newcommand{\fho}[1]{\vphantom{\mid}\setbox0\hbox{00}\hbox to \wd0{\hss\ensuremath{#1}\hss}}
\newcommand{\wt}{\widetilde}
\newcommand{\wh}{\widehat}
\newcommand{\ad}[1]{\mathrm{ad}_{#1}}
\newcommand{\tr}{\mathrm{tr}\,}
\newcommand{\GL}{\mathrm{GL}}
\newcommand{\SL}{\mathrm{SL}}
\newcommand{\SO}{\mathrm{SO}}
\newcommand{\Sp}{\mathrm{Sp}}
\newcommand{\Mat}{\mathrm{Mat}}
\newcommand{\Stab}{\mathrm{Stab}}

\newcommand{\vfi}{\varphi}
\newcommand{\teta}{\vartheta}
\newcommand{\Bfi}{\Phi}
\newcommand{\Fp}{\mathbb{F}}
\newcommand{\Rp}{\mathbb{R}}
\newcommand{\Zp}{\mathbb{Z}}
\newcommand{\Cp}{\mathbb{C}}
\newcommand{\ut}{\mathfrak{u}}
\newcommand{\at}{\mathfrak{a}}
\newcommand{\nt}{\mathfrak{n}}
\newcommand{\mt}{\mathfrak{m}}
\newcommand{\htt}{\mathfrak{h}}
\newcommand{\spt}{\mathfrak{sp}}
\newcommand{\rt}{\mathfrak{r}}
\newcommand{\rad}{\mathfrak{rad}}
\newcommand{\bt}{\mathfrak{b}}
\newcommand{\gt}{\mathfrak{g}}
\newcommand{\vt}{\mathfrak{v}}
\newcommand{\pt}{\mathfrak{p}}
\newcommand{\Xt}{\mathfrak{X}}
\newcommand{\Po}{\mathcal{P}}
\newcommand{\Uo}{\EuScript{U}}
\newcommand{\Fo}{\EuScript{F}}
\newcommand{\Do}{\EuScript{D}}
\newcommand{\Eo}{\EuScript{E}}
\newcommand{\Iu}{\mathcal{I}}
\newcommand{\Mo}{\mathcal{M}}
\newcommand{\Nu}{\mathcal{N}}
\newcommand{\Ro}{\mathcal{R}}
\newcommand{\Co}{\mathcal{C}}
\newcommand{\Lo}{\mathcal{L}}
\newcommand{\Ou}{\mathcal{O}}
\newcommand{\Uu}{\mathcal{U}}
\newcommand{\Au}{\mathcal{A}}
\newcommand{\Vu}{\mathcal{V}}
\newcommand{\Bu}{\mathcal{B}}
\newcommand{\Sy}{\mathcal{Z}}
\newcommand{\Sb}{\mathcal{F}}
\newcommand{\Gr}{\mathcal{G}}
\newcommand{\rtc}[1]{C_{#1}^{\mathrm{red}}}

\author{Mikhail V. Ignatyev}

\date{}
\title{\Large{Gradedness of the set of rook placements in $A_{n-1}$}\mbox{$\vphantom{1}$}\footnotetext{
The work on Section~\ref{sect:Kerop_map} was performed at the NRU HSE with the financial support from the Russian Science Foundation, grant no. 16--41--01013.}} \maketitle

\begin{center}
\begin{tabular}{p{15cm}}
\small{\textsc{Abstract}. A rook placement is a subset of a root system consisting of positive roots with pairwise non-positive inner products. To each rook placement in a root system one can assign the coadjoint orbit of the Borel subgroup of a reductive algebraic group with this root system. Degenerations of such orbits induce a natural partial order on the set of rook placements. We study combinatorial structure of the set of rook placements in $A_{n-1}$ with respect to a slightly different order and prove that this poset is graded.}
\end{tabular}
\end{center}

\sect{Introduction}

\sst Denote by $G=\GL_n(\Cp)$ the group of all invertible $n\times n$ matrices over the complex numbers. Let $B$ be the Borel subgroup of $G$ consisting of all invertible upper-triangular matrices, $U$ be the unipotent radical of $B$ (it consists of all upper-triangular matrices with $1$'s on the diagonal), and $T$ be the subgroup of all invertible diagonal matrices (it is the maximal torus of $G$ contained in $B$). Next, let $\bt$ and $\nt$ be the Lie algebras of $B$ and $U$ respectively.

Let $\Phi$ be the root system of $G$ with respect to~$T$, $\Phi^+$ be the set of positive roots with respect to~$B$, $\Delta$ be the set of simple roots, and $W$ be the Weyl group of $\Phi$ (for basic facts on algebraic groups and root systems, see \cite{Bourbaki}, \cite{Humphreys} and \cite{Humpreys2}). The root system $\Phi$ is of type $A_{n-1}$; as usual, we identify the set of positive roots with the subset of the Euclidean space $\Rp^n$ of the form $$A_{n-1}^+=\{\epsi_i-\epsi_j,~1\lee i<j\lee n\}.$$ Under this identification, $\Delta$ consists of the roots $\alpha_i=\epsi_i-\epsi_{i+1}$, $1\lee i\lee n-1$ ($\{\epsi_i\}_{i=1}^n$ is the standard basis of $\Rp^n$).

\defi{A \emph{rook placement} is a subset $D\subseteq\Phi^+$ such that $(\alpha,\beta)\lee0$ for all distinct $\alpha$, $\beta\in D$. (Here $(\cdot,\cdot)$ denotes the standard inner product on $\Rp^n$.)}

\exam{Let $n=6$. Below we draw the rook placement $D=\{\epsi_1-\epsi_3,~\epsi_2-\epsi_6,\epsi_3-\epsi_5\}$. If a root $\epsi_i-\epsi_j$ is contained in $D$ then we put the symbol $\otimes$ in the $(j,i)$th entry of the $n\times n$ chessboard. If we interpret these symbols as rooks then it follows from the definition that the rooks do not hit each other.
\begin{equation*}\predisplaypenalty=0
\mymatrix{
\Bot{2pt}\pho& \pho& \pho& \pho& \pho& \pho\\
\pho& \Lft{2pt}\Bot{2pt}\pho& \pho& \pho& \pho& \pho\\
\otimes& \pho& \Lft{2pt}\Bot{2pt}\pho& \pho& \pho& \pho\\
\pho& \pho& \pho& \Lft{2pt}\Bot{2pt}\pho& \pho& \pho\\
\pho& \pho& \otimes& \pho& \Lft{2pt}\Bot{2pt}\pho& \pho\\
\pho& \otimes& \pho& \pho& \pho& \Lft{2pt}\pho\\
}\end{equation*}}

We denote the set of all rook placement in $A_{n-1}$ by $\Ro(n)$. Further, let $\Iu(n)$ be the set of all orthogonal rook placements. Below we describe two closely related partial orders on these sets.

The Lie algebra $\nt$ has the basis $\{e_{\alpha},~\alpha\in\Phi^+\}$ consisting of the root vectors: for $\alpha=\epsi_i-\epsi_j$, $e_{\alpha}$ is nothing but the elementary matrix $e_{i,j}$. Denote by $\{e_{\alpha}^*,~\alpha\in\Phi^+\}$ the dual basis of the dual space $\nt^*$. Given a rook placement $D$, put $$f_D=\sum_{\beta\in D}e_{\beta}^*\in\nt^*.$$ The group $B$ acts on its Lie algebra $\bt$ by the adjoint action, and $\nt$ is an invariant subspace. Hence one has the dual action of the groups $B$ and $U$ on the space $\nt^*$; we call this action \emph{coadjoint}. We say that the $B$-orbit $\Omega_D\subset\nt^*$ of the linear form $f_D$ is \emph{associated} with the rook placement $D$.

Such orbits play an important role in the A.A. Kirillov's orbit method \cite{Kirillov1}, \cite{Kirillov2} describing representations of $B$ and $U$. For $D\in\Iu(n)$, such orbits were studied by A.N. Panov in \cite{Panov1} and by me in \cite{Ignatyev1}. One can define analogues of such orbits for other root systems, see \cite{Ignatyev2}, \cite{Ignatyev3}, \cite{Ignatyev4} for the case of $\Iu(n)$. For arbitrary rook placements in $\Ro(n)$, such orbits were considered in \cite{IgnatyevVasyukhin1}; see also \cite{Andre1}, \cite{AndreNeto1}, where C. Andre and A. Neto used rook placements to construct so-called supercharacter theory for the group $U$. Note that in \cite{Melnikov1}, \cite{Melnikov2}, A. Melnikov studied the adjoint $B$-orbits of elements of the form $\sum_{\beta\in D}e_{\beta}$, $D\in\Iu(n)$.

Given a subset $A\subseteq\nt^*$, we will denote by $\overline{A}$ its closure with respect to the Zarisski topology. There exists a natural partial order on the set $\Ro(n)$ induced by the degenerations of associated orbits: we will write $D_1\lee_B D_2$ if $\Omega_{D_1}\subseteq\overline{\Omega}_{D_2}$. We need to introduce one more partial order on the set of rook placements. Namely, given an arbitrary $D\in\Ro(n)$, denote by $R_D$ the $n\times n$ matrix defined by
\begin{equation*}
(R_D)_{i,j}=\begin{cases}\#\{\epsi_a-\epsi_b\in D\mid a\lee j,~b\gee i\},&\text{if }i>j,\\
0&\text{otherwise}.
\end{cases}
\end{equation*}
Put $D_1\lee D_2$ if $(R_{D_1})_{i,j}\lee(R_{D_2})_{i,j}$ for all $i$, $j$.

\exam{Let $n=4$, $D_1=\{\epsi_1-\epsi_2,~\epsi_2-\epsi_4\}$, $D_2=\{\epsi_1-\epsi_3,~\epsi_2-\epsi_4\}$. Then
\begin{equation*}
D_1=\mymatrix{
\Bot{2pt}\pho& \pho& \pho& \pho\\
\otimes& \Lft{2pt}\Bot{2pt}\pho& \pho& \pho\\
\pho& \pho& \Lft{2pt}\Bot{2pt}\pho& \pho\\
\pho& \otimes& \pho& \Lft{2pt}\pho\\
}~,~R_{D_1}=\begin{pmatrix}
0&0&0&0\\
1&0&0&0\\
0&1&0&0\\
0&1&1&0
\end{pmatrix},~
D_2=\mymatrix{
\Bot{2pt}\pho& \pho& \pho& \pho\\
\pho& \Lft{2pt}\Bot{2pt}\pho& \pho& \pho\\
\otimes& \pho& \Lft{2pt}\Bot{2pt}\pho& \pho\\
\pho& \otimes& \pho& \Lft{2pt}\pho\\
}~,~R_{D_2}=\begin{pmatrix}
0&0&0&0\\
1&0&0&0\\
1&2&0&0\\
0&1&1&0
\end{pmatrix}.
\end{equation*}
We conclude that $D_1\lee D_2$. On the other hand, it is easy to check that $D_1\nleqslant_B D_2$, see\break \cite[Re\-mark~1.6 (iii)]{IgnatyevVasyukhin1}, so these two partial orders on $\Ro(n)$ do not coincide.}

Nevertheless, it turns out that these orders are closely related to each other. Precisely, given rook placements $D_1$, $D_2\in\Ro(n)$, it follows from $D_1\lee_B D_2$ that $D_1\lee D_2$ \cite[Theo\-rem~1.5]{IgnatyevVasyukhin1}. Furthermore, if $D_1$, $D_2\in\Iu(n)$ then the conditions $D_1\lee_B D_2$ and $D_1\lee D_2$ are equivalent \cite[Theo\-rem~1.7]{Ignatyev1}. Besides, given a rook placement $$D=\{\epsi_{i_1}-\epsi_{j_1},~\ldots,~\epsi_{i_l}-\epsi_{j_l}\},$$ we denote by $w_D\in S_n$ the permutation, which is equal to the product of transpositions $$w_D=(i_1,j_1)\ldots(i_l,j_l).$$ Now, both of the conditions above (for orthogonal rook placements $D_1$, $D_2$) are equivalent to the condition that $w_{D_1}$ is less or equal to $w_{D_2}$ with respect to the Bruhat order \cite[Theorem 1.1]{Ignatyev1}. Similar facts are true for orthogonal rook placements in the root system $C_n$, see \cite{Ignatyev2}. Note that these results are in some sense ``dual'' to A. Melnikov's results.

In the paper \cite{Incitti2}, F. Incitti studied the order on $\Iu(n)$ induced by the Bruhat order on the elements $w_D$, $D\in\Iu(n)$, from purely combinatorial point of view (see also \cite{Incitti1} for other classical root systems). In particular, given an orthogonal rook placement $D$, he explicitly described the set of its immediate predecessors (it consists of $D'\in\Iu(n)$ such that there exists an edge from $D'$ to $D$ in the Hasse diagram of this poset). The set of immediate predecessors for the partial order $\lee$ on $\Iu(n)$ and $\Ro(n)$ was described by me in \cite[Lemmas 3.6, 3.7, 3.8]{Ignatyev1} and by A.S. Vasyukhin and me in \cite[Theo\-rem~3.3]{IgnatyevVasyukhin1} respectively. (In the case of $\Iu(n)$, the set of immediate predecessors for $\lee$ coincides with the set described by F.~Incitti, which implies that those two partial orders coincide.)

Furthermore, F. Incitti proved that the poset $\Iu(n)$ is graded and calculated its M$\mathrm{\ddot o}$bius function. Recall that a finite poset $X$ is called \emph{graded} if it has the greatest and the lowest elements and all maximal chains in $X$ have the same length. Gradedness is equivalent to the existence of a rank function. By definition, it is a (unique) function $\rho$ on $X$, which value on the lowest element is zero, such that if $x$ is an immediate predecessor of $y$ then $\rho(y)=\rho(x)+1$ In \cite[Theorem 5.2]{Incitti2}, F. Incitti constructed the rank function on $\Iu(n)$. As the main result of this paper, we prove the gradedness of the poset $\Ro(n)$.

The main tool used in the proof is so-called Kerov placements (see \cite{Kerov1}). To each rook placement $D\in\Ro(n)$ one can assign a certain orthogonal rook placement $K(D)\in\Iu(2n-2)$. We prove that if rook placements $D_1$ is an immediate predecessor of $D_2$ in $\Ro(n)$ then $K(D_1)$ is an immediate predecessor of $K(D_2)$ in $\Iu(2n-2)$ (and vice versa), see Theorem~\ref{theo:A_n_Kerov}. As a corollary, we construct a rank function on $\Ro(n)$ and prove the gradedness of this poset, see Corollary~\ref{coro:A_n_graded}.

The structure of the paper is as follows. In the next section we describe the set of immediate predecessors of a given rook placement for $\Iu(n)$ and $\Ro(n)$. In the third section we introduce the Kerov map $$K\colon\Ro(n)\to\Iu(2n-2)$$ and show that it preserves the property ``to be an immediate predecessor''. This allows us to use F. Incitti's results to construct a rank function on $\Ro(n)$, which imlies the gradedness of this poset.

\medskip\textsc{Acknowledgements}. A part of this work (Section~\ref{sect:imm_pre}) was done during my stay at the Oberwolfach Research Institute for Mathematics in February--March 2018 via the program ``Research in pairs'' together with Alexey Petukhov. I thank Alexey for fruitful discussions. I also thank the Oberwolfach Research Institute for Mathematics for hospitality and financial support.

\sect{Immediate predecessors}\label{sect:imm_pre}

To prove that the set $\Ro(n)$ is graded with respect to the partial order introduced above, we need to describe the set of immediate predecessors of a given rook placement in $\Ro(n)$ and $\Iu(n)$. Such a description for $\Ro(n)$ was provided in \cite{IgnatyevVasyukhin1}, while for $\Iu(n)$ it was presented in F. Incitti's work \cite{Incitti2}. Recall that a rook placement $D\in\Ro(n)$ is called an \emph{immediate predecessor} of a rook placement $T\in\Ro(n)$ if $D<T$ and there are no $S\in\Ro(n)$ such that $D<S<T$. (As usual,  $D<T$ means that $D\lee T$ and $D\neq T$.) In other words, there exists an oriented edge from $D$ to $T$ in the Hasse diagram of the poset $\Ro(n)$. The definition of immediate predecessors for $\Iu(n)$ is literally the same.

We denote the set of all immediate predecessors in $\Ro(n)$ (respectively, in $\Iu(n)$) of a rook placement $D\in\Ro(n)$ (respectively, of an orthogonal rook placement $D\in\Iu(n)$) by $L_{\Ro}(D)$ (respectively, by $L_{\Iu}(D)$). This set consists of rook placements of several types, which we will describe now. First, we will consider the set $L_{\Ro}(D)$ in details.

It is convenient to introduce the following notation. We will write simply $(i,j)$ instead of $\epsi_j-\epsi_i$, $i>j$. Besides, for each $k$ from 1 to $n$, we put
\begin{equation*}
\Ro_k=\{(k,s)\in\Phi^+\mid1\leq s<k\},~\Co_k=\{(r,k)\in\Phi^+\mid {j<k\leq n}\}.
\end{equation*}

\defi{The sets $\Ro_k$, $\Co_k$ are called the $k$th \emph{row} and the $k$th \emph{column} of $\Phi^+$ respectively. We will write $\row(\alpha)=k$ and $\col(\alpha)=k$ if $\alpha\in\Ro_k$ and $\alpha\in\Co_k$ respectively. Note that, for $D\in\Ro(n)$, one has
\begin{equation*}\label{formula:Ro_Co_RP}
|D\cap\Ro_k|\leq1\text{ and }|D\cap\Co_k|\leq1\text{ for all }1\leq k\leq n.
\end{equation*}
Furthermore, if $D\in\Iu(n)$ then
\begin{equation*}\label{formula:Ro_Co_inv}
|D\cap(\Ro_k\cup\Co_k)|\leq1\text{ for all }1\leq k\leq n.
\end{equation*}}

There exists a natural partial order on the set of positive roots: given $\alpha,~\beta\in\Phi^+$, by definition, $\alpha\leq\beta$ if $\beta-\alpha$ is a (probably, empty) sum of positive roots. In the other words, $$(a,b)\leq(c,d)\text{ if $c\geq a$ and $d\leq b.$}$$ Given a rook placement~$D\in\Ro(n)$, denote by $\wt M(D)$ the set of minimal roots from $D$ (with respect to~$\lee$). Now, we set
\begin{equation*}
\begin{split}
M_{\Ro}(D)&=\{(i,j)\in\wt M(D)\mid D\cap\Ro_k\neq\varnothing\text{ and }D\cap\Co_k\neq\varnothing\text{ for all }j<k<i\},\\
N_{\Ro}^-(D)&=\{D_{(i,j)}^-,(i,j)\in M_{\Ro}(D)\},
\end{split}
\end{equation*}
where $D_{(i,j)}^-=D\setminus\{(i,j)\}$.

Next, fix a root $(i,j)\in D$. Denote $$m=\min\{k\mid j<k<i\text{ и }D\cap\Co_k=\varnothing\}.$$
Suppose that such a number $m$ exists. Assume that $D\cap\Ro_k\neq\varnothing$ for all $k$ from $j+1$ to $m$. Assume, in addition, that there are no $(p,q)\in D$ such that $(i,j)>(p,q)$ and $(i,m)\not>(p,q)$. The set of all roots $(i,j)\in D$ satisfying these conditions is denoted by $A^{\Ro}_{\to}(D)$; given $(i,j)\in A^{\Ro}_{\to}(D)$, we put $$D_{(i,j)}^{\to,\Ro}=(D\setminus\{(i,j)\})\cup\{(i,m)\}.$$
Similarly, suppose that there exists a number $$m=\max\{k\mid j<k<i\text{ и }D\cap\Ro_k=\varnothing\}.$$ Assume also that $D\cap\Co_k\neq\varnothing$ for $m+1\leq k\leq i-1$ and that there are no $(p,q)\in D$ such that $(i,j)>(p,q)$ и $(m,j)\not>(p,q)$. Denote the set of all such $(i,j)$'s by $A^{\Ro}_{\uparrow}$; given $(i,j)\in A^{\Ro}_{\uparrow}$, we put $$D_{(i,j)}^{\uparrow,\Ro}=(D\setminus\{(i,j)\})\cup\{(m,j)\}.$$

Now, let $B^{\Ro}_{(i,j)}(D)$ be the set of roots $(\alpha,\beta)\in D$ such that $(\alpha,\beta)>(i,j)$ and there are no $(p,q)\in D$ satisfying $(i,j)<(p,q)<(\alpha,\beta)$. For each $(\alpha,\beta)\in B^{\Ro}_{(i,j)}(D)$ we set $$D_{(i,j)}^{(\alpha,\beta),\Ro}=(D\setminus\{(i,j),(\alpha,\beta)\})\cup\{(i,\beta),(\alpha,j)\}.$$
By definition, let $$N^0_{\Ro}(D)=\left\{D_{(i,j)}^{\uparrow,\Ro},~(i,j)\in A^{\Ro}_{\uparrow}\right\}\cup\left\{D_{(i,j)}^{\to,\Ro},~(i,j)\in A^{\Ro}_{\to}\right\}\cup\bigcup_{(i,j)\in D}\left\{D_{(i,j)}^{(\alpha,\beta),\Ro},(\alpha,\beta)\in B^{\Ro}_{(i,j)}(D)\right\}.$$

\exam{Let $n=8$ and $D=\{(3,1),(6,2),(7,3),(5,4),(8,5)\}$. Clearly, $M_{\Ro}(D)=\{(5,4)\}$, $(8,5)\in A^{\Ro}_{\to}$, $(3,1)\in A^{\Ro}_{\uparrow}$ and $(6,2)\in B^{\Ro}_{(5,4)}(D)$. On the picture below we draw the rook placements $D$, $D_{(5,4)}^{(6,2),\Ro}$, $D_{(3,1)}^{\uparrow,\Ro}$ and $D_{(8,5)}^{\to,\Ro}$.
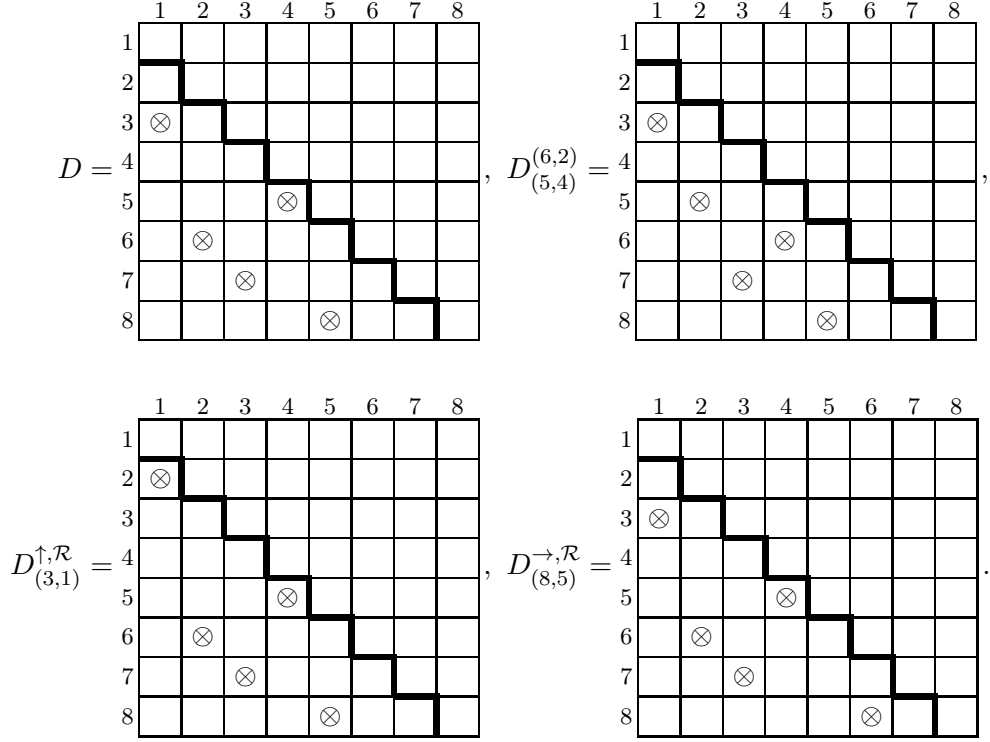
\begin{equation*}
\begin{split}
D&=\mymatrix{ \pho& \pho& \pho& \pho& \pho& \pho& \pho& \pho\\
\Top{2pt}\Rt{2pt}\pho& \pho& \pho& \pho& \pho& \pho& \pho& \pho\\
\otimes& \Top{2pt}\Rt{2pt}\pho& \pho& \pho& \pho& \pho& \pho& \pho\\
\pho& \pho& \Top{2pt}\Rt{2pt}\pho& \pho& \pho& \pho& \pho& \pho\\
\pho& \pho& \pho& \Top{2pt}\Rt{2pt}\otimes& \pho& \pho& \pho& \pho\\
\pho& \otimes& \pho& \pho& \Top{2pt}\Rt{2pt}\pho& \pho& \pho& \pho\\
\pho& \pho& \otimes& \pho& \pho& \Top{2pt}\Rt{2pt}\pho& \pho& \pho\\
\pho& \pho& \pho& \pho& \otimes& \pho& \Top{2pt}\Rt{2pt}\pho& \pho\\}\,,~
D_{(5,4)}^{(6,2)}=\mymatrix{ \pho& \pho& \pho& \pho& \pho& \pho& \pho& \pho\\
\Top{2pt}\Rt{2pt}\pho& \pho& \pho& \pho& \pho& \pho& \pho& \pho\\
\otimes& \Top{2pt}\Rt{2pt}\pho& \pho& \pho& \pho& \pho& \pho& \pho\\
\pho& \pho& \Top{2pt}\Rt{2pt}\pho& \pho& \pho& \pho& \pho& \pho\\
\pho& \otimes& \pho& \Top{2pt}\Rt{2pt}\pho& \pho& \pho& \pho& \pho\\
\pho& \pho& \pho& \otimes& \Top{2pt}\Rt{2pt}\pho& \pho& \pho& \pho\\
\pho& \pho& \otimes& \pho& \pho& \Top{2pt}\Rt{2pt}\pho& \pho& \pho\\
\pho& \pho& \pho& \pho& \otimes& \pho& \Top{2pt}\Rt{2pt}\pho& \pho\\}\,,\\\\
D_{(3,1)}^{\uparrow,\Ro}&=\mymatrix{ \pho& \pho& \pho& \pho& \pho& \pho& \pho& \pho\\
\Top{2pt}\Rt{2pt}\otimes& \pho& \pho& \pho& \pho& \pho& \pho& \pho\\
\pho& \Top{2pt}\Rt{2pt}\pho& \pho& \pho& \pho& \pho& \pho& \pho\\
\pho& \pho& \Top{2pt}\Rt{2pt}\pho& \pho& \pho& \pho& \pho& \pho\\
\pho& \pho& \pho& \Top{2pt}\Rt{2pt}\otimes& \pho& \pho& \pho& \pho\\
\pho& \otimes& \pho& \pho& \Top{2pt}\Rt{2pt}\pho& \pho& \pho& \pho\\
\pho& \pho& \otimes& \pho& \pho& \Top{2pt}\Rt{2pt}\pho& \pho& \pho\\
\pho& \pho& \pho& \pho& \otimes& \pho& \Top{2pt}\Rt{2pt}\pho& \pho\\}\,,~
D_{(8,5)}^{\to,\Ro}=\mymatrix{ \pho& \pho& \pho& \pho& \pho& \pho& \pho& \pho\\
\Top{2pt}\Rt{2pt}\pho& \pho& \pho& \pho& \pho& \pho& \pho& \pho\\
\otimes& \Top{2pt}\Rt{2pt}\pho& \pho& \pho& \pho& \pho& \pho& \pho\\
\pho& \pho& \Top{2pt}\Rt{2pt}\pho& \pho& \pho& \pho& \pho& \pho\\
\pho& \pho& \pho& \Top{2pt}\Rt{2pt}\otimes& \pho& \pho& \pho& \pho\\
\pho& \otimes& \pho& \pho& \Top{2pt}\Rt{2pt}\pho& \pho& \pho& \pho\\
\pho& \pho& \otimes& \pho& \pho& \Top{2pt}\Rt{2pt}\pho& \pho& \pho\\
\pho& \pho& \pho& \pho& \pho& \otimes& \Top{2pt}\Rt{2pt}\pho& \pho\\}\,.
\end{split}
\end{equation*}}

Next, fix a root $(i,j)\in D$, and consider a pair $(\alpha,\beta)\in\Zp\times\Zp$. Suppose that $i>\beta\geq\alpha>j$, $D\cap\Ro_{\alpha}=D\cap\Co_{\beta}=\varnothing$, $D\cap\Ro_k\neq\varnothing$, $D\cap\Co_k\neq\varnothing$ for all $\alpha<k<\beta$, and the conditions $(p,q)\in D$, $(i,j)>(p,q)$, $(\alpha,j)\not>(p,q)$ imply $(i,\beta)>(p,q)$. Moreover, assume that if $\alpha\neq\beta$ then $D\cap\Ro_{\beta}\neq\varnothing$ and $D\cap\Co_{\alpha}\neq\varnothing$. Denote the set of all such pairs $(\alpha,\beta)$ by $C^{\Ro}_{(i,j)}(D)$. For an arbitrary pair $(\alpha,\beta)\in C^{\Ro}_{(i,j)}(D)$, we put $$D_{(i,j)}^{\alpha,\beta,\Ro}=(D\setminus\{(i,j)\})\cup\{(i,\beta),(\alpha,j)\}.$$
By definition, let $$N^+_{\Ro}(D)=\bigcup_{(i,j)\in D}\left\{D_{(i,j)}^{\alpha,\beta,\Ro},(\alpha,\beta)\in C^{\Ro}_{(i,j)}(D)\right\}.$$

\exam{Let $n=6$ and $D=\{(4,1),(6,2),(5,4)\}$, then $(3,3)\in C^{\Ro}_{(6,2)}(D)$. On the picture below we draw the rook placements $D$ and $D_{(6,2)}^{3,3,\Ro}$.
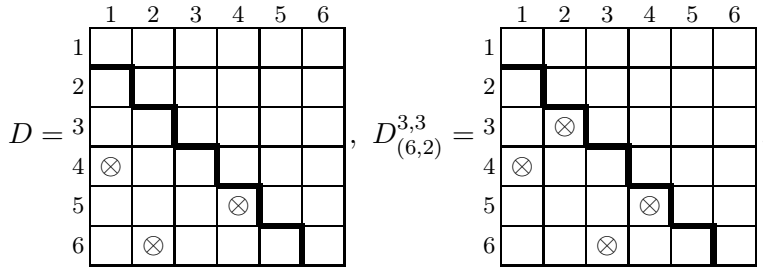
\begin{equation*}
D=\mymatrix{ \pho& \pho& \pho& \pho& \pho& \pho\\
\Top{2pt}\Rt{2pt}\pho& \pho& \pho& \pho& \pho& \pho\\
\pho& \Top{2pt}\Rt{2pt}\pho& \pho& \pho& \pho& \pho\\
\otimes& \pho& \Top{2pt}\Rt{2pt}\pho& \pho& \pho& \pho\\
\pho& \pho& \pho& \Top{2pt}\Rt{2pt}\otimes& \pho& \pho\\
\pho& \otimes& \pho& \pho& \Top{2pt}\Rt{2pt}\pho& \pho\\}\,,~
D_{(6,2)}^{3,3}=\mymatrix{ \pho& \pho& \pho& \pho& \pho& \pho\\
\Top{2pt}\Rt{2pt}\pho& \pho& \pho& \pho& \pho& \pho\\
\pho& \Top{2pt}\Rt{2pt}\otimes& \pho& \pho& \pho& \pho\\
\otimes& \pho& \Top{2pt}\Rt{2pt}\pho& \pho& \pho& \pho\\
\pho& \pho& \pho& \Top{2pt}\Rt{2pt}\otimes& \pho& \pho\\
\pho& \pho& \otimes& \pho& \Top{2pt}\Rt{2pt}\pho& \pho\\}\,.
\end{equation*}}

Finally, we set
\begin{equation*}
N_{\Ro}(D)=N_{\Ro}^-(D)\cup N^0_{\Ro}(D)\cup N^+_{\Ro}(D).
\end{equation*}
The set of immediate predecessors of a given rook placement from $\Ro(n)$ is described as follows.

\theoc{\cite[Theorem 3.3]{IgnatyevVasyukhin1}}{Let $D\in\Ro(n)$.\label{theo:cov_rel} Then $L_{\Ro}(D)=N(D)$.}

Now we turn to the description of immediate predecessors for $\Iu(n)$. Given an orthogonal rook placement~$D\in\Ro(n)$, put
\begin{equation*}
\begin{split}
M_{\Iu}(D)&=\{(i,j)\in\wt M(D)\mid D\cap(\Ro_k\cup\Co_k)\neq\varnothing\text{ for all }j<k<i\},\\
N_{\Iu}^-(D)&=\{D_{(i,j)}^-,(i,j)\in M_{\Iu}(D)\},
\end{split}
\end{equation*}
where $D_{(i,j)}^-=D\setminus\{(i,j)\}$, as above.

Let $D\in\Iu(n)$, $(i,j)\in D$. Denote $$m=\min\{k\mid j<k<i\text{ и }D\cap\Co_k=D\cap\Ro_k=\varnothing\}.$$ Suppose that such a number $m$ exists. Assume that there are no $(p,q)\in D$ such that $(i,j)>(p,q)$ and $(i,m)\not>(p,q)$. The set of all $(i,j)\in D$ satisfying these conditions is denoted by $A^{\Iu}_{\to}(D)$; given $(i,j)\in A^{\Iu}_{\to}(D)$, we set $$D_{(i,j)}^{\to,\Iu}=(D\setminus\{(i,j)\})\cup\{(i,m)\}.$$
Similarly, suppose that there exists $$m=\max\{k\mid j<k<i\text{ и }D\cap\Ro_k=D\cap\Co_k=\varnothing\},$$ and there are no $(p,q)\in D$ such that $(i,j)>(p,q)$ and $(m,j)\not>(p,q)$. The set of all such $(i,j)$'s is denoted by $A^{\Iu}_{\uparrow}$; given $(i,j)\in A^{\Iu}_{\uparrow}$, we set $$D_{(i,j)}^{\uparrow,\Iu}=(D\setminus\{(i,j)\})\cup\{(m,j)\}.$$

Next, let $B^{\Iu}_{(i,j)}(D)$ be the set of roots $(\alpha,\beta)\in D$ such that $j<\beta<i<\alpha$, $$D\cap(\Ro_r\cup\Co_r)\neq\varnothing$$ for all $\beta<r<i$ and there are no $(p,q)\in D$ for which $j<q<\beta<p<i$ or $\beta<q<i<p<\alpha$ (in other words, for which $(i,j)>(p,q)$ and $(\beta,j)\ngtr(p,q)$, or $(\alpha,\beta)>(p,q)$ and $(\alpha,i)\ngtr(p,q)$). To each $(\alpha,\beta)\in B^{\Iu}_{(i,j)}(D)$ we assign the set
\begin{equation*}
D_{(i,j)}^{(\alpha,\beta),\Iu}=
(D\setminus\{(i,j),(\alpha,\beta)\})\cup\{(\beta,j),(\alpha,i)\}.
\end{equation*}
Now, let
\begin{equation*}
\begin{split}
N^0_{\Iu}(D)&=\left\{D_{(i,j)}^{\uparrow,\Iu},~(i,j)\in A^{\Iu}_{\uparrow}\right\}\cup\left\{D_{(i,j)}^{\to,\Iu},~(i,j)\in A^{\Iu}_{\to}\right\}\\
&\cup\bigcup_{(i,j)\in D}\left\{D_{(i,j)}^{(\alpha,\beta),\Ro},(\alpha,\beta)\in B^{\Ro}_{(i,j)}(D)\right\}\cup\bigcup_{(i,j)\in D}\left\{D_{(i,j)}^{(\alpha,\beta),\Iu},(\alpha,\beta)\in B^{\Iu}_{(i,j)}(D)\right\}.
\end{split}
\end{equation*}

\exam{If $n=8$, $D=\{(5,1),(6,2),(8,4)\}$, then $(8,4)\in B_{6,2}^{\Iu}(D)$, hence
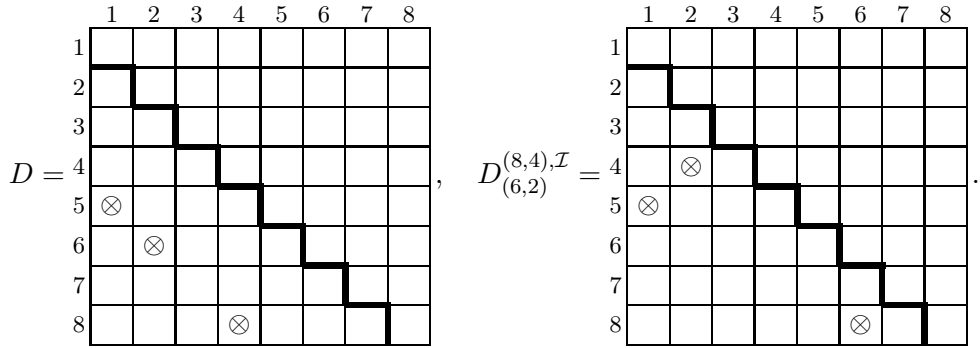
\begin{equation*}D=
\mymatrix{ \pho& \pho& \pho& \pho& \pho& \pho& \pho& \pho\\
\Top{2pt}\Rt{2pt}\pho& \pho& \pho& \pho& \pho& \pho& \pho& \pho\\
\pho& \Top{2pt}\Rt{2pt}\pho& \pho& \pho& \pho& \pho& \pho& \pho\\
\pho& \pho& \Top{2pt}\Rt{2pt}\pho& \pho& \pho& \pho& \pho& \pho\\
\otimes& \pho& \pho& \Top{2pt}\Rt{2pt}\pho& \pho& \pho& \pho& \pho\\
\pho& \otimes& \pho& \pho& \Top{2pt}\Rt{2pt}\pho& \pho& \pho& \pho\\
\pho& \pho& \pho& \pho& \pho& \Top{2pt}\Rt{2pt}\pho& \pho& \pho\\
\pho& \pho& \pho& \otimes& \pho& \pho& \Top{2pt}\Rt{2pt}\pho& \pho\\
}\,,\quad D_{(6,2)}^{(8,4),\Iu}=
\mymatrix{ \pho& \pho& \pho& \pho& \pho& \pho& \pho& \pho\\
\Top{2pt}\Rt{2pt}\pho& \pho& \pho& \pho& \pho& \pho& \pho& \pho\\
\pho& \Top{2pt}\Rt{2pt}\pho& \pho& \pho& \pho& \pho& \pho& \pho\\
\pho& \otimes& \Top{2pt}\Rt{2pt}\pho& \pho& \pho& \pho& \pho& \pho\\
\otimes& \pho& \pho& \Top{2pt}\Rt{2pt}\pho& \pho& \pho& \pho& \pho\\
\pho& \pho& \pho& \pho& \Top{2pt}\Rt{2pt}\pho& \pho& \pho& \pho\\
\pho& \pho& \pho& \pho& \pho& \Top{2pt}\Rt{2pt}\pho& \pho& \pho\\
\pho& \pho& \pho& \pho& \pho& \otimes& \Top{2pt}\Rt{2pt}\pho& \pho\\
}\,.\end{equation*}}

Besides, denote by $C^{\Iu}_{i,j}(D)$ the set of pairs $(\alpha,\beta)\in{\mathbb{Z}}\times{\mathbb{Z}}$ such that $i>\beta>\alpha>j$, $$D\cap(\Ro_{\alpha}\cup\Co_{\alpha})=D\cap(\Ro_{\beta}\cup\Co_{\beta})=\varnothing,$$ $D\cap(\Ro_k\cup\Co_k)\neq\varnothing$ for all $\beta>k>\alpha$, and if $(p,q)\in D$, $(i,j)>(p,q)$, $(\alpha,j)\ngtr(p,q)$ then $(i,\beta)>(p,q)$. For each pair $(i,j)\in C^{\Iu}_{(i,j)}(D)$, we put
\begin{equation*}
D_{(i,j)}^{\alpha,\beta,\Iu}=
(D\setminus\{(i,j)\})\cup\{(i,\beta),(\alpha,j)\}.
\end{equation*}

\exam{Let $n=8$, $D=\{(4,1),(8,2),(7,6)\}$, then $(3,5)\in
C^{\Iu}_{(8,2)}(D)$, so
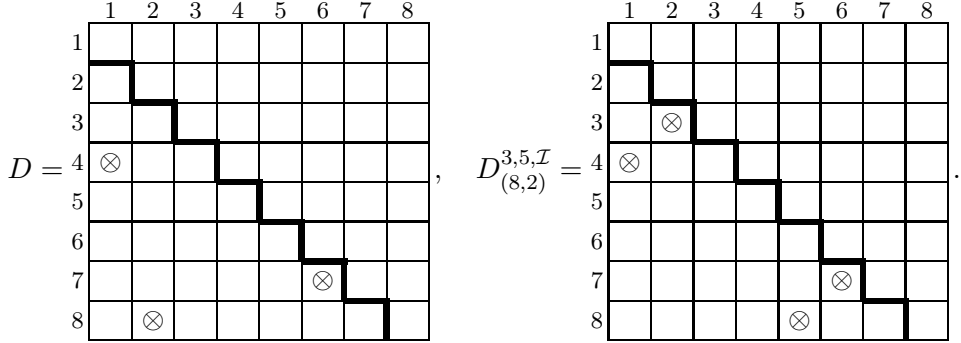
\begin{equation*}D=
\mymatrix{ \pho& \pho& \pho& \pho& \pho& \pho& \pho& \pho\\
\Top{2pt}\Rt{2pt}\pho& \pho& \pho& \pho& \pho& \pho& \pho& \pho\\
\pho& \Top{2pt}\Rt{2pt}\pho& \pho& \pho& \pho& \pho& \pho& \pho\\
\otimes& \pho& \Top{2pt}\Rt{2pt}\pho& \pho& \pho& \pho& \pho& \pho\\
\pho& \pho& \pho& \Top{2pt}\Rt{2pt}\pho& \pho& \pho& \pho& \pho\\
\pho& \pho& \pho& \pho& \Top{2pt}\Rt{2pt}\pho& \pho& \pho& \pho\\
\pho& \pho& \pho& \pho& \pho& \Top{2pt}\Rt{2pt}\otimes& \pho& \pho\\
\pho& \otimes& \pho& \pho& \pho& \pho& \Top{2pt}\Rt{2pt}\pho& \pho\\
}\,,\quad D_{(8,2)}^{3,5,\Iu}=
\mymatrix{ \pho& \pho& \pho& \pho& \pho& \pho& \pho& \pho\\
\Top{2pt}\Rt{2pt}\pho& \pho& \pho& \pho& \pho& \pho& \pho& \pho\\
\pho& \Top{2pt}\Rt{2pt}\otimes& \pho& \pho& \pho& \pho& \pho& \pho\\
\otimes& \pho& \Top{2pt}\Rt{2pt}\pho& \pho& \pho& \pho& \pho& \pho\\
\pho& \pho& \pho& \Top{2pt}\Rt{2pt}\pho& \pho& \pho& \pho& \pho\\
\pho& \pho& \pho& \pho& \Top{2pt}\Rt{2pt}\pho& \pho& \pho& \pho\\
\pho& \pho& \pho& \pho& \pho& \Top{2pt}\Rt{2pt}\otimes& \pho& \pho\\
\pho& \pho& \pho& \pho& \otimes& \pho& \Top{2pt}\Rt{2pt}\pho& \pho\\
}\,.\end{equation*}}
Finally, we denote
\begin{equation*}
\begin{split}
N^+_{\Iu}(D)&=\bigcup_{(i,j)\in D}\left\{D_{(i,j)}^{\alpha,\beta,\Iu},(\alpha,\beta)\in C^{\Iu}_{(i,j)}(D)\right\},\\
N_{\Iu}(D)&=N_{\Ro}^-(D)\cup N^0_{\Iu}(D)\cup N^+_{\Iu}(D).
\end{split}
\end{equation*}

Immediate predecessors in $\Iu(n)$ are described by the following F. Incitti's theorem (see also \cite[Subsection 2.4]{Ignatyev1}).

\theoc{\cite[Theorem 5.1]{Incitti2}}{Let $D\in\Iu(n)$.\label{theo:cov_rel_Iu} Then $L_{\Iu}(D)=N_{\Iu}(D)$.}

\sect{Kerov map and the main result}\label{sect:Kerop_map}

In this section, we introduce our main technical tool, Kerov orthogonal rook placements, and, using them, prove that $\Ro(n)$ is graded.

\defi{Let $n\geq3$, and $D$ be a rook placement from $\Ro(n)$. A \emph{Kerov rook placement} corresponding to $D$ is, by definition, the orthogonal rook placement $K(D)\in\Iu(2n-2)$ construdcted by the following rule: if $$D=\{(i_1,j_1),\ldots,(i_s,j_s)\},$$ then $$K(D)=(2i_1-2,2j_1-1)\ldots(2i_r-2,2j_r-1).$$ (Kerov rook placements were introduced in the paper \cite{Kerov1}). We call the map $K\colon\Ro(n)\to\Iu(2n-2)$ given by the rule $D\mapsto K(D)$ the \emph{Kerov map}.}

\exam{If $n=8$ and $D=\{(3,1),(6,2),(7,3),(5,4),(8,6)\}\in\Ro(8)$, then
\begin{equation*}
\begin{split}
K(D)&=(4,1)\cdot(10,3)\cdot(12,5)\cdot(8,7)\cdot(14,11)\\
&=\begin{pmatrix}
1&2&3&4&5&6&7&8&9&10&11&12&13&14\\
4&2&10&1&12&6&8&7&9&3&14&5&13&11
\end{pmatrix}\in\Iu(14).
\end{split}
\end{equation*}}

The following proposition is evident.

\mprop{Let $D$, $T\in\Ro(n)$. Then the conditions $T\leq D$ and $K(T)\leq K(D)$ are equivalent.}

The following theorem plays the crucial role in the proof of the main result.

\theop{Let $D,~T\in\Ro(n)$ be rook placements\label{theo:A_n_Kerov}. Then the conditions $T\in L_{\Ro}(D)$ and\break $K(T)\in L_{\Iu}(K(D))$ are equivalent.}{Clearly, $K(T)\in L_{\Iu}(D)$ implies $T\in L_{\Ro}(D)$. Indeed, since there are no orthogonal involutions from $\Iu(2n-2)$ between $K(T)$ and $K(D)$, we conclude that, in particular, there are no Kerov involutions between them. It remains to prove that the converse is also true.

Assume that $T\in L_{\Ro}(D)$. By Theorem~\ref{theo:cov_rel}, this is equivalent to $T\in N_{\Ro}(D)=N^-(D)\cup N_{\Ro}^0(D)\cup N_{\Ro}^+(D)$. We will consider these variants case-by-case.

First, suppose that $T\in N_{\Ro}^-(D)$. This means that $T=D_{(i,j)}^-$ for a certain root $(i,j)\in M(D)$. Automatically, $K(T)=K(D)\setminus\{(2i-2,2j-1)\}$. It follows immediately from $(i,j)\in\wt M(D)$ that $(2i-2,2j-1)\in\wt M(K(D))$. Since $(i,j)\in M(D)$, we see that $D\cap\Ro_k$ and $D\cap\Co_k$ are nonempty if $i<k<j$. This shows that $K(D\cap\Ro_{2k-2})$ and $K(D)\cap\Co_{2k-1}$ are nonempty for all such $k$. Thus, $(2i-2,2j-1)\in M(K(D))$, i.e., $K(T)\in N_{\Iu}^-(K(D))$. By Theorem~\ref{theo:cov_rel_Iu}, $K(T)\in L_{\Iu}(K(D))$.

Next, assume that $T\in N_{\Ro}^0(D)$. If $T=D_{(i,j)}^{(\alpha,\beta),\Ro}$ for some $(i,j)\in D$, $(\alpha,\beta)\in\Bu_{(i,j)}^{\Ro}(D)$, then it is easy to see that $(2\alpha-2,2\beta-1)\in\Bu_{(2i-2,2j-1)}^{\Ro}(K(D))$ и $K(T)=K(D)_{(2i-2,2j-1)}^{(2\alpha-2,2\beta-1),\Ro}\in N_{\Ro}^0(K(D))$, hence $K(T)\in N_{\Iu}^0(D)\subset L_{\Iu}(K(D))$.
Now consider the case when $T=D_{(i,j)}^{\to,\Ro}$ for some $(i,j)\in A_{\to}^{\Ro}$. (The case $T=D_{(i,j)}^{\uparrow,\Ro}$, $(i,j)\in A_{\uparrow}^{\Ro}$ can be considered similarly.) Let $T=(D\setminus\{(i,j)\})\cup\{(i,m)\}$, then $K(T)=(K(D)\setminus\{(2i-2,2j-1)\})\cup\{(2i-2,2m-1)\}$. Since there are no root in $D$ which is less than $(i,j)$ but not less than $(i,m)$, we have a similar condition for $K(D)$. Since $D\cap\Co_k\neq\varnothing$ при $j<k<m$, one has $K(D)\cap\Co_{2k-1}\neq\varnothing$ for such $k$. On the other hand, $D\cap\Ro_k$ is nonempty for $j<k\leq m$, so $K(D)\cap\Ro_{2k-2}$ is also nonempty for such $k$. Thus, $K(D)\cap(\Ro_k\cup\Co_k)\neq\varnothing$ for $2j-1<k<2m-1$, which means that $(2i-2,2j-1)\in A_{\to}^{\Iu}$ and $K(T)=K(D)_{(2i-2,2j-1)}^{\to,\Iu}$. Hence, by Theorem~\ref{theo:cov_rel_Iu}, $K(T)\in L_{\Iu}(K(D))$, as required.

Finally, suppose that $T\in N_{\Ro}^+(D)$, i.e., $T=D_{(i,j)}^{\alpha,\beta,\Ro}$ for certain $(i,j)\in D$ and $(\alpha,\beta)\in C_{(i,j)}^{\Ro}(D)$. Since $i>\beta\geq\alpha>j$, we have $2i-2>2\beta-1>2\alpha-2>2j-1$. It follows from $D\cap\Ro_{\alpha}=D\cap\Co_{\beta}=\varnothing$ that $K(D)\cap\Ro_{2\alpha-2}=K(D)\cap\Co_{2\beta-1}=\varnothing$. Since $K(D)$ is a Kerov rook placement, the condition $K(D)\cap\Co_{2\alpha-2}=K(D)\cap\Ro_{2\beta-1}=\varnothing$ is satisfied automatically. If $\alpha=\beta$ then there is nothing to prove. If $\beta>\alpha$ then $D\cap\Ro_k\neq\varnothing$ and $D\cap\Co_k\neq\varnothing$ for all $k$ from $\alpha+1$ to $\beta-1$, hence $K(D)\cap\Ro_{2k-2}\neq\varnothing$ and $K(D)\cap\Co_{2k-1}\neq\varnothing$ for all such $k$. Furthermore, $D\cap\Ro_{\beta}$ and $D\cap\Co_{\alpha}$ are nonempty, which implies that $K(D)\cap\Ro_{2\beta-2}$ and $=D\cap\Co_{2\alpha-1}$ are also nonempty. Thus, we obtain $K(D)\cap(\Ro_k\cap\Co_k)\neq\varnothing$ для всех $k$ от $2\alpha-1$ до $2\beta-2$, как и требовалось. Таким образом, $(2\alpha-2,2\beta-1)\in C^{\Iu}_{(2i-2,2j-1)}(D)$ и $K(T)=K(D)_{(2i-2,2j-1)}^{2\alpha-2,2\beta-1,\Iu}$. Theorem~\ref{theo:cov_rel_Iu} guarantees that $K(T)\in L_{\Iu}(K(D))$. The proof is complete.}

\corop{For\label{coro:A_n_graded} each $n\gee2$ the poset $\Ro(n)$ is graded with the rank function
\begin{equation*}
\rho(D) = \dfrac{l(w_{K(D)})+|D|}{2},
\end{equation*}
where $l(w)$ is the length of a permutation $w$ in the corresponding symmetric group.
}{As we mentioned in the introduction, F. Incitii showed that the set $\Iu(2n-2)$ of orthogonal rook placements id graded. Precisely \cite[Theorem 5.3.2]{Incitti1}, the rank function on this poset has the form
\begin{equation*}
\rho(D)=\dfrac{l(w_D)+|D|}{2}.
\end{equation*}
Applying Theorem~\ref{theo:A_n_Kerov}, we see that the restriction of this rank function to $K(\Ro(n))$ in fact provided the rank function of the required form on $\Ro(n)$. This concludes the proof.}

\bigskip\textsc{Mikhail V. Ignatyev}\par\textsc{Samara National Research University, ul. Ak. Pavlova, 1, 443011, Samara, Russia}

\emph{E-mail}: \texttt{mihail.ignatev@gmail.com}

\end{document}